\documentclass[10pt,a4paper]{article}
\usepackage{amsmath}
\usepackage{amssymb}
\usepackage{amsbsy}
\usepackage{amsfonts}
\usepackage{graphicx}
\usepackage{sectsty}
\usepackage{verbatim}
\usepackage{enumitem}
\usepackage{hyperref}
\usepackage{ifthen}
\usepackage{tikz}

\RequirePackage{ifthen}
\makeatletter
\newcommand{\logmessage}[1]{\@latex@warning{#1}}
\makeatother
\IfFileExists{../Bibliographie/Users_Guide.txt}{
  }{
  \IfFileExists{../../Bibliographie/Users_Guide.txt}{
    }{
    \IfFileExists{../../../Bibliographie/Users_Guide.txt}{
      }{
        \IfFileExists{../../../../Bibliographie/Users_Guide.txt}{
          }{
            \IfFileExists{../../../../../Bibliographie/Users_Guide.txt}{
              }{
        \logmessage{Directory 'Bibliographie' not found}
          }}}}}
\numberwithin{equation}{section}
\newtheorem{remark}{Remark}[section]

\newcommand{\field}[1]{\ensuremath{\mathbb{#1}}}
\newcommand{\R}{\field{R}}

\renewcommand{\i}{\mathrm i}

\newcommand{\F}{{\cal F}}

\newcommand{\intr}{\int_{-\infty}^\infty}

\newcommand{\abs}[1]{\lvert#1\rvert}

\makeatletter
\newcommand{\ignore}{\logmessage{Text ignored}\@gobble}
\makeatother

\newcommand{\smalf}{\par\smallskip\noindent}
\newcommand{\medlf}{\par\medskip\noindent}

\newcommand{\cmplx}{\mathbb{C}}

\title{Simultaneous Reconstructions of Absorption Density and Wave Speed with
Photoacoustic Measurements
\thanks{The work has been supported by the Austrian Science Fund (FWF)
within the national research network Photo\-acoustic Imaging in Biology and
Medicine, project S10505-N20.}}

\author{Andreas Kirsch \thanks{Karlsruhe Institute of Technology,
Department of Mathematics, Kaiserstra\ss{}e 89, 76128 Karlsruhe, Germany}
\and Otmar Scherzer \thanks{Computational Science Center, University
of Vienna, Nordbergstra\ss{}e 15, A-1090 Vienna, Austria, and
Johann Radon Institute for Computational and Applied Mathematics (RICAM),
Austrian Academy of Sciences, Altenbergerstra\ss{}e 69,
A-4040 Linz, Austria ({\tt otmar.scherzer@univie.ac.at}).}
}

\begin{document}
\maketitle
\begin{abstract}
In this paper we propose an approach for \emph{simultaneous} identification of
the \emph{absorption density} and the \emph{speed of sound} by photoacoustic
measurements.
Experimentally our approach can be realized with sliced photoacoustic experiments.
The mathematical model for such an experiment is developed and exact reconstruction
formulas for both parameters are presented.
\end{abstract}
\section{Introduction}

\begin{figure}[h]
\centering\includegraphics[width=0.70\textwidth]{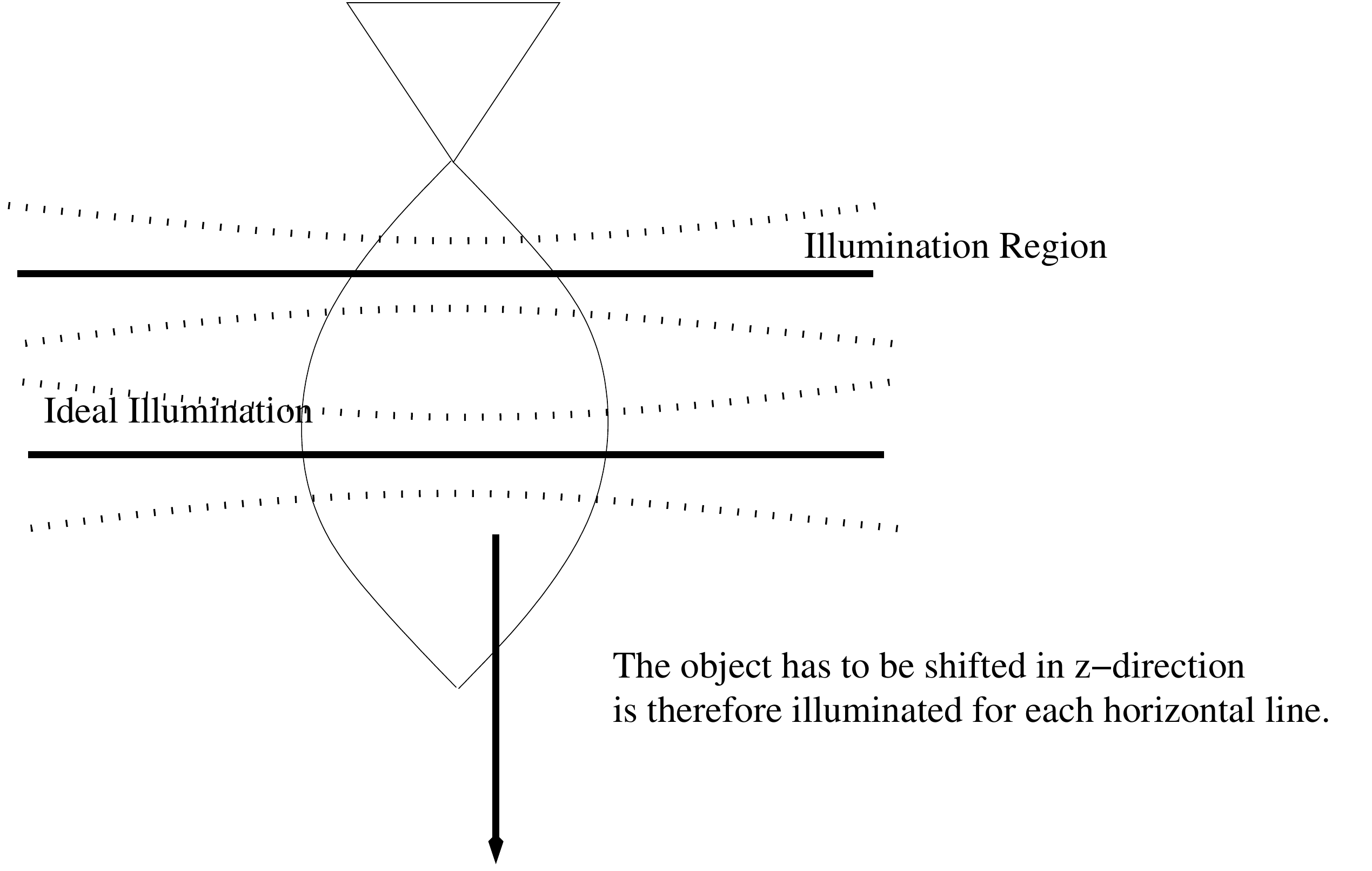}
\caption{\label{fig:fish} Conventional sectional photoacoustic imaging.}
\end{figure}

\begin{figure}[h]
\centering\includegraphics[width=0.70\textwidth]{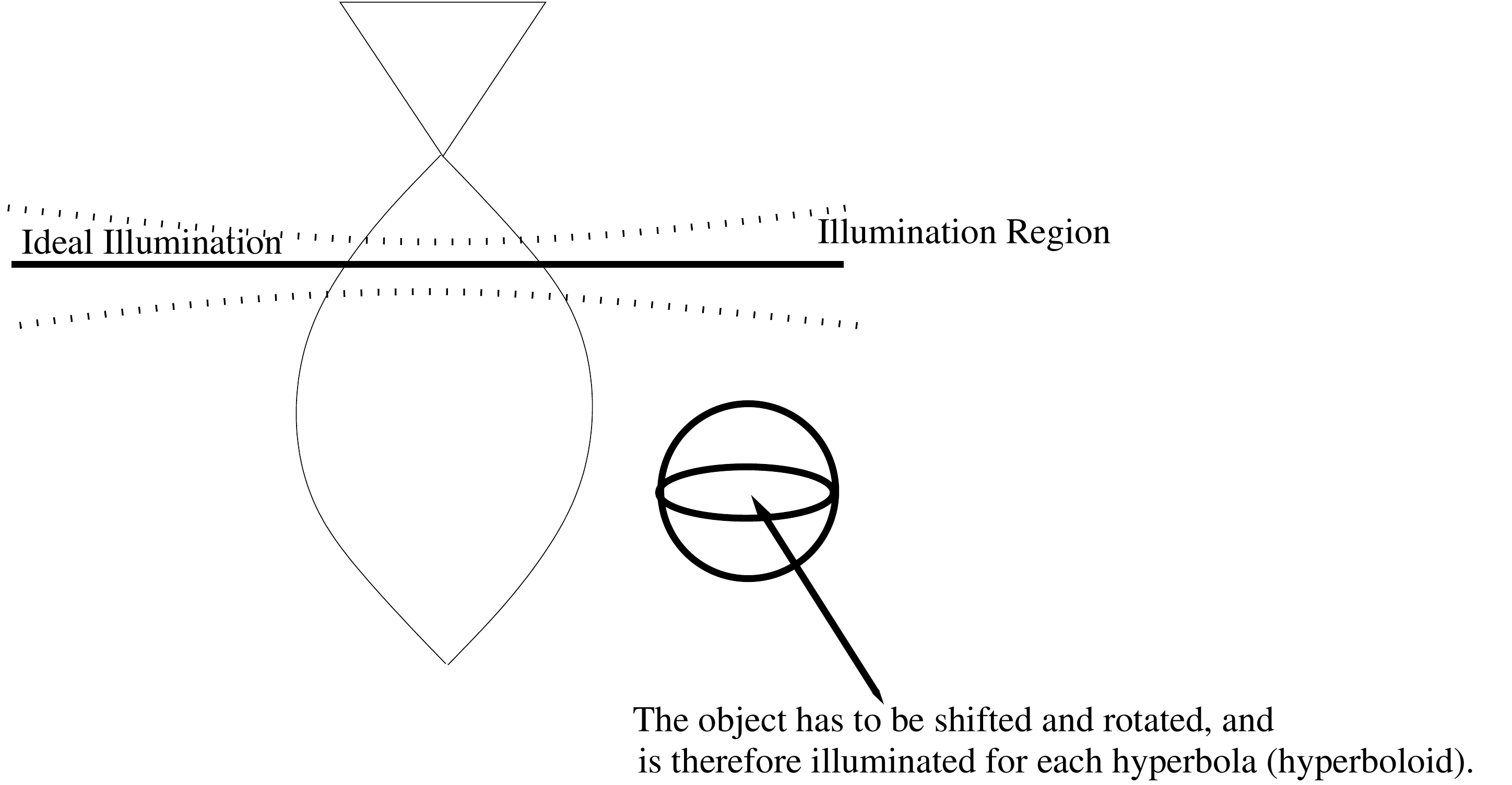}
\caption{\label{fig:fish_two} Sectional imaging in all directions produces enough data to reconstruct both imaging functions.}
\end{figure}

In this paper we propose an approach for \emph{simultaneous} identification of
the \emph{absorption density} and the \emph{speed of sound} by photoacoustic
measurements.
\smalf
In standard photoacoustic experiments the object is \emph{uniformly} illuminated
by a short electromagnetic impulse. Recently \emph{sectional} photoacoustic
imaging techniques \cite{MaTarNtzRaz09,RazDieNtz09,GratPasNusPal11,GratPasNusPal11_2} 
and according reconstructions techniques \cite{ElbSchSchu11_report} have been
developed. In sectional experiments thin hyperbola shaped like regions of the specimen 
are imaged (see Figure \ref{fig:fish}) sequentially. 
Experimentally, one can perform sectional photoacoustic imaging by focused illumination 
combined with focusing detectors. The focused illumination is achieved by using cylindrical 
lenses in front of the object. Moreover, contemporary focusing ultrasonic detectors have 
a spherical or cylindrical shape, thus the detector surface plays the role of the acoustic
lens. The generated ultrasonic wave is refracted by a suitable acoustic lens such that 
out-of-plane (center of the hyperbola shaped regions) signals are generally weak and 
can be neglected. Thus essentially only signals emerging from the imaging plane are 
collected at the detector. This justifies that the illumination can be assumed restricted 
to a single plane (line in 2D). However, such a model requires a low scattering coefficient 
of the sample (which model organism specimens like the Zebra fish have).
Opposed to \emph{sectional 3D photoacoustic imaging}, where stacks of
complementary two-dimensional projection images are produced, the proposed
approach for simultaneous imaging consists in performing overlapping sliced
imaging by rotation and translation of the specimen. This, also generates enough
data for reconstructing the two independent parameter functions, speed of sound
and absorption density, (see Figure \ref{fig:fish_two}).
By now the approach presented here is by far from being experimentally economical
in the sense that a quick dimension analysis shows that we acquire much too many
data - in $\R^3$ we reconstruct two 3-dimensional functions from eventually six
dimensional data. It is the goal of this paper, however, not to present the most
economical approach, but to show that it is possible to derive exact
reconstruction formulas for both imaging parameter functions.
\smalf
The literature on reconstruction formulas and back-projection algorithms for
photoacoustic imaging is vast. Wang et al.\ developed reconstruction formulas
for cylindrical, spherical, and planar measurement geometries in a series of
papers~\cite{XuFenWan02,XuWan02a,XuWan02b,XuXuWan02}, and recently many more
algorithms based on reconstruction formulas have been developed
(see the survey~\cite{KucKun08}). As it becomes transparent below we can
make use of reconstruction formulas from photoacoustics--essentially
we make extensive use of inversion formulas for the spherical mean operator. 
Exact reconstruction formulas for the wavespeed function in ultrasound reflectivity 
tomography are based on the \emph{Born approximation} to the wave equation, which is 
valid for moderately varying speed of sound. Reconstruction formulas for ultrasound 
reflectivity tomography have already been derived by Norton \& Linzer \cite{Nor80,NorLin81} 
and are also based on inversion formulas for the spherical mean operator. The possibility of 
exact inversion in both fields supports to derive exact inversion formulas for both parameters.
\smalf
The proposed method is a hybrid and quantitative imaging methods (see \cite{BalUhl10,Kuc11,Bal11}) 
for some recent surveys. Most closely related to our approach is the work of Stefanov and Uhlmann \cite{SteUhl11a} 
which presented a photoacoustic experiment for recovering either the sound speed or absorption density. 
A substitute to our work is also the recovery of the absorption density function for inhomogeneous wave speed
\cite{AgrKuc07,HriKucNgu08,SteUhl09}.
\smalf
The reconstruction formulas for simultaneous imaging utilize techniques from
reflectivity imaging and photoacoustic imaging. In the current state of research
we can provide exact reconstruction formulas, but when we try to be economical,
i.e., using less slicing experiments, we would require reconstruction formulas, which
have not been developed so far. This is highlighted in Section \ref{sec:model}.
The outline of this paper is as follows: In Section \ref{sec:model} we discuss
the mathematical model and discuss the non-economicality of proposed model,
which leaves room for further improvement of the results. In Section
\ref{sec:Fourier} we transform the model into the Fourier domain, from which in
Sections \ref{sec:2D}, \ref{sec:3D} exact reconstruction formulas are derived in
2D, 3D, respectively.
The Appendix provides the exact
definitions used in this paper.

\section{Model}
\label{sec:model}

In $\R^n$, we consider the following Cauchy problem for the wave equation:
\begin{eqnarray*}
\frac{1}{c^2}\,\partial_{tt}\tilde{u}\ -\ \Delta\tilde{u} & = & 0\quad
\text{in }\R^n\times\R_{>0}\,, \\
\tilde{u}(x,0) & = & f(x)\,\delta_{r,\theta}(x)\quad\text{in }\R^n\,, \\
\partial_t\tilde{u}(x,0) & = & 0\quad\text{in }\R^n\,.
\end{eqnarray*}
where $c=c(x)$ denotes the speed of sound and $\delta_{r,\theta}(x)=
\delta(\text{dist}(x,E(r,\theta)))$ where $E(r,\theta)$ is the
$(n-1)$-dimensional hyperplane with distance $r$ from the origin and orientation
$\theta\in S^{n-1}$. This is the photoacoustic equation with sliced illumination
in the plane $E(r,\theta)$. It is the ultimate aim to reconstruct $c$ and $f$
from measurements of $\tilde{u}$ on some hypersurface $\Gamma$ (see below).
\smalf
We consider the Born approximation; that is, we expand $\tilde{u}$ formally with
respect to the contrast function  $q:=1/c^2-1$ and consider only terms of order
at most $q$. This leads to the decomposition $\tilde{u}\approx u+v$ where
$u=u^{r,\theta}$ is the solution of the wave equation
\begin{eqnarray*}
\partial_{tt} u\ -\ \Delta u & = & 0\quad\text{in }\R^n\times\R_{>0}\,,
\nonumber \\
u(x,0) & = & f(x)\,\delta_{r,\theta}(x)\quad\text{in }\R^n\,,\\
\partial_t u(x,0) & = & 0\quad\text{in }\R^n\,, \nonumber
\end{eqnarray*}
and $v=v^{r,\theta}$ solves
\begin{eqnarray*}
\partial_{tt} v -\ \Delta v & = & -q(x)\,\partial_{tt} u\quad
\mbox{in }\R^n\times\R_{>0}\,, \nonumber \\
v(x,0) & = & 0\quad\text{in }\R^n\,, \\
\partial_t v(x,0) & = & 0\quad\text{in }\R^n\,. \nonumber
\end{eqnarray*}
The dependence of $v$ on $r,\theta$ is through the source function $u$.
\smalf
It is our goal to reconstruct $q$ and $f$ from measurements of
\begin{equation} \label{eq:measurements}
m^{r,\theta}(x,t)\ =\ u^{r,\theta}(x,t)\ +\ v^{r,\theta}(x,t)\,,\quad
(x,t)\in\Gamma\times (0,T)\,,
\end{equation}
where $\Gamma$ is a $(n-1)$-dimensional hypersurface in $\R^n$. Currently
everything is fixed to complete measurements on the \emph{whole} surface $\Gamma$.
\smalf
We make the following assumptions:
\begin{enumerate}
\item $\Gamma=\partial B$ for some open and connected set $B\subset\R^n$ with
smooth boundary. As typically in photoacoustics, we consider $B$ to be a sphere,
a circle, or a halfplane.
\item $f$ is the sum of some known initial distribution $f_0$ with compact
support and some unknown term $f_1$ to be determined.
\item The supports of $f_1$ and $q$ are both contained in $\Omega$ for some
bounded domain $\Omega\subset B$,
\item $f_0,f_1$ and $q$ are smooth and $f(x)=f_0(x)+f_1(x)\neq 0$ for all
$x\in\overline{B}$.
\end{enumerate}

\begin{figure}[h]
\centering\includegraphics[width=0.50\textwidth]{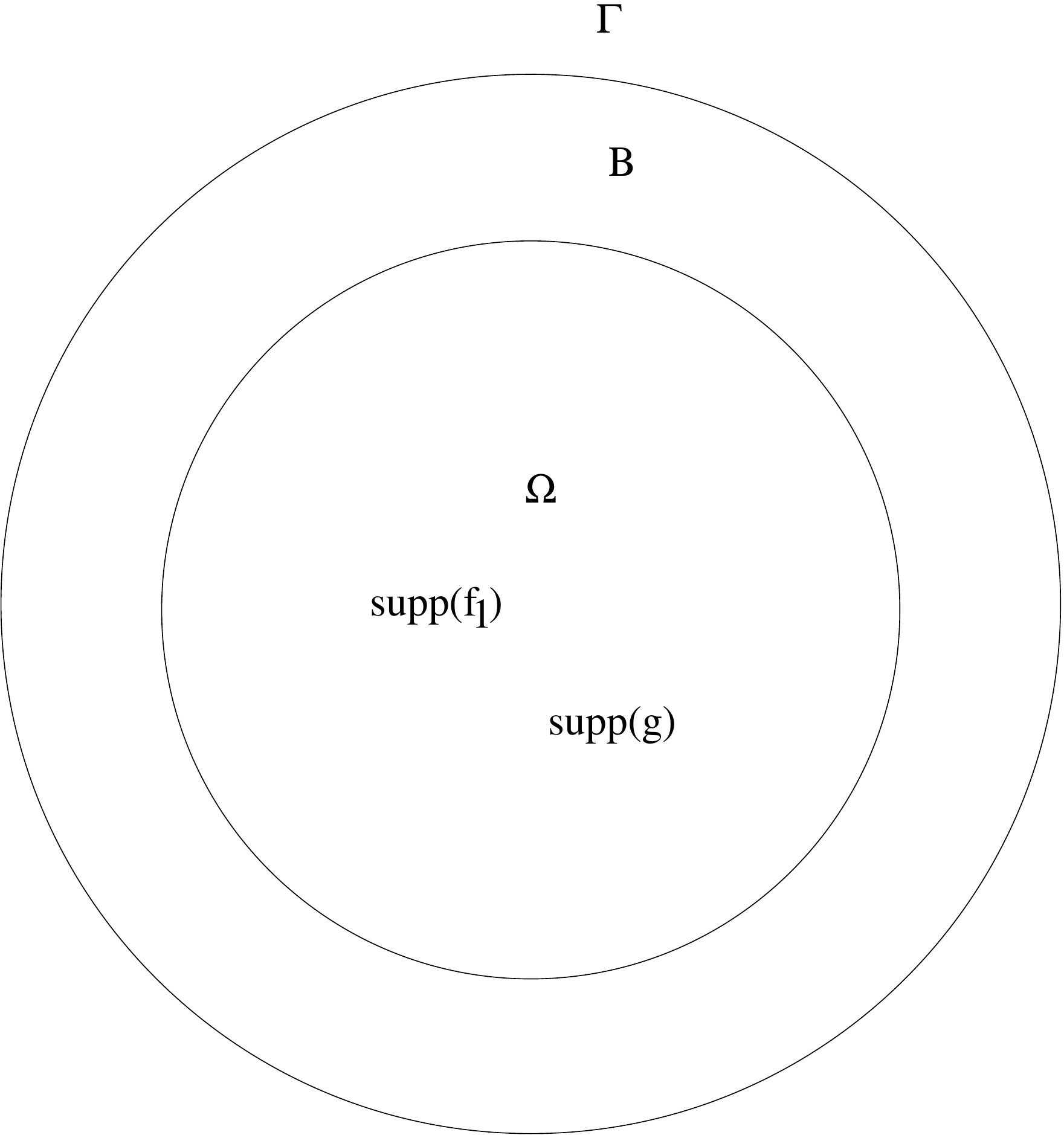}
\caption{\label{fig1} Schematic representation of the domains and supports. }
\end{figure}

\subsection*{Dimensionality Analysis}

In the $n$-dimensional setting, we record measurements (cf.
(\ref{eq:measurements})) for every $(x,t)\in\Gamma\times(0,T)$ and every
$(r,\theta)\in(0,\infty)\times S^{n-1}$. That is, the recorded data are
$2n$-dimensional. The data to be recovered are $f_1$ and $q$, which
are two $n$-dimensional functions. In general we think that we record too many
data for the purpose of reconstructing $f_1$ and $q$ -- however, since we rely
on Radon transforms techniques for exact inversion, we accept this disadvantage
for the mathematical studies.

\section{Fourier Reconstruction Formulas}
\label{sec:Fourier}

In this section we derive exact reconstruction formulas for $f_1$ and $g$. First
we derive a general inversion formula, which is then evaluated differently in 2D
and 3D.
\smalf
In the following we omit the superscripts $(r,\theta)$ for the sake of
convenience of notation. Let
$$ \hat{u}(x,k)\ =\ \F u(x,k)\ =\ \int_0^\infty u(x,t)\,e^{\i k t}\,dt\,,\quad
x\in\R^n\,,\ k\in\R\,, $$
be the Fourier transform - note that $u$ vanishes for $t <0$. 
We assume that $u$ is bounded in $\R^n\times\R_{>0}$
such that $u(x,\cdot)\in L^2(0,\infty)$ for every $x\in\R^n$. Then we note that
$\hat{u}$ has the following properties (by the theorem of Paley-Wiener):
\begin{itemize}
\item[(a)] $u(x,\cdot)$ has a holomorphic extension into $\cmplx_+:=\{z\in\cmplx:
\Im z>0\}$ for all $x\in\R^n$,
\item[(b)] for every $x\in\R^n$ it holds that $\int_{-\infty}^\infty
\bigl|\hat{u}(x,k_1+ik_2)-\hat{u}(x,k_1)\bigr|^2dk_1\longrightarrow 0$ as $k_2$
tends to zero,
\item[(c)] $\hat{u}(\cdot,k)$ is bounded in $\R^n$ for all $k\in\cmplx_+$.
\end{itemize}
We have
\begin{eqnarray*}
\hat{u}(x,k) & = & \F u(x,k)\ =\
\int_0^\infty u(x,t)\,e^{\i k t}\,dt \\
& = & \frac{1}{\i k}\left[\left. u(x,t)\,e^{\i kt}\right|_{t=0}^\infty\
-\ \int_0^\infty \partial_t u(x,t)\,e^{\i kt}dt\right] \\
& = & -\frac{1}{\i k}\,u(x,0)\ -\ \frac{1}{(\i k)^2}
\left[\left.\partial_tu(x,t)\,e^{\i kt}\right|_{t=0}^\infty\ -\
\int_0^\infty\partial_{tt}u(x,t)\,e^{\i k t}dt\right] \\
& = & -\frac{1}{\i k}\,u(x,0)\ -\ \frac{1}{k^2}\,\Delta\hat{u}(x,k)\,,
\end{eqnarray*}
or in other words
\begin{equation} \label{eq:Helmholtz}
\Delta\hat{u}(x,k)\ +\ k^2\hat{u}(x,k)\ =\ \i k\,u(x,0)\,.
\end{equation}
Moreover, we have
\begin{eqnarray*}
\hat{v}(x,k) & = & \F v(x,k)\ =\
-\frac{1}{k^2}\int_0^\infty v_{tt}(x,t)\,e^{\i kt}\,dt \\
& = & -\frac{1}{k^2}\int_0^\infty\left[\Delta v(x,t)-
q(x)\,\partial_{tt}u(x,t)\right]\,dt \\
& = & -\frac{1}{k^2}\,\Delta\hat{v}(x,k)\ +\ \frac{q(x)}{k^2}
\int_0^\infty\partial_{tt}u(x,t)\,e^{\i kt}dt\\
& = & -\frac{1}{k^2}\,\Delta\hat{v}(x,k)\ +\ \frac{q(x)}{k^2}\,
\Delta\hat{u}(x,k)\,,
\end{eqnarray*}
and thus
\begin{eqnarray*}
\Delta\hat{v}(x,k)\ +\ k^2\hat{v}(x,k)\ & = & q(x)\,\Delta\hat{u}(x,k) \\
& = & -k^2q(x)\,\hat{u}(x,k)\ +\ \i k\,q(x)\,u(x,0)\,.
\end{eqnarray*}
Let $\Phi_k$ be the (radiating) fundamental solution of the Helmholtz equation
$\Delta\hat{u}+k^2\hat{u}=0$ in $\R^n$ for $k\in\cmplx$ with $\Im z\geq 0$; that is,
in particular for $n=2,3$,
$$ \Phi_k(x,y)\ :=\ \left\{\begin{array}{cl}
\displaystyle \frac{\exp(ik|x-y|}{4\pi|x-y|} & \mbox{for }n=3\,, \\[3mm]
\displaystyle \frac{\i}{4}\,H_0^{(1)}(k|x-y|) & \mbox{for }n=2\,,
\end{array}\right.\quad x\not=y\,, $$
where $H_0^{(1)}$ denotes the Hankel function of the first kind and order zero.
Then \emph{one particular} solution of (\ref{eq:Helmholtz}) is given by
\begin{eqnarray*}
\hat{u}(x,k) & = & -\i k\int_{\R^n}u(y,0)\,\Phi_k(x,y)\,dy \\
& = & -\i k\int_{y\in E(r,\theta)}\!\!\!f(y)\,\Phi_k(x,y)\,ds(y)\,.
\end{eqnarray*}
We note that this particular solution satisfies the properties (a), (b), (c) from
the beginning of this section. Any other solution is of the form
$\hat{u}(\cdot,k)+w(\cdot,k)$ where $w$ satisfies $\Delta w+k^2w=0$ in $\R^n$.
The requirement that the solution satisfies the properties (a), (b), (c) from
the beginning of this section implies that $w$ vanishes. Indeed, for $k\in
\cmplx_+$ the function $w(\cdot,k)$ has to be bounded by property (c), therefore
its Fourier transform (in the distributional sense) with respect to $x\in\R^n$
satisfies
$$ \mathcal{F}_xw(y,k)\,\bigl[k^2-|y|^2\bigr]\ =\ 0\mbox{ for all }y\in\R^n $$
which implies $\mathcal{F}_xw(\cdot,k)=0$ because $k^2-|y|^2$ does not vanish for
$k\in\cmplx_+$. Therefore, also $w(\cdot,k)=0$ in $\R^n$ for all $k\in\cmplx_+$
and, by the continuity property (b), also $w(\cdot,k)=0$ in $\R^n$ for all
$k\in\cmplx$ with $\Im\geq 0$.
\smalf
We also have
\begin{eqnarray*}
& & \hat{v}(x,k) \\
& = & k^2\int_{\R^n}q(y)\,\hat{u}(y,t)\,\Phi_k(x,y)\,dy\ -\
    \i k\int_{\R^n}q(y)\,f(y)\,\delta_{r,\theta}(y)\,\Phi_k(x,y)\,dy \\
& = & -\i k^3\int_{\R^n}q(y)\int_{z\in E(r,\theta)}\!\!\!f(z)\,\Phi_k(y,z)\,
ds(z)\,\Phi_k(x,y)\,dy \\
& & \qquad -\ \i k\int_{z\in E(r,\theta)}\!\!\!q(z)\,f(z)\,\Phi_k(y,z)\,ds(z)\\
& = & -\i k\int\limits_{z\in E(r,\theta)}\!\!\!\!\!\!\!\!\! f(z)\left[
k^2\int_{\R^n} q(y)\,\Phi_k(y,z)\,\Phi_k(x,y)\,dy\ +\
q(z)\,\Phi_k(x,z)\right]ds(z)\,.
\end{eqnarray*}
In summary, we have
\begin{eqnarray*}
\hat m^{r,\theta}(x,k) & = & \hat{u}^{r,\theta}(x,k)\ +\ \hat{v}^{r,\theta}(x,k)
\ =\ \hat{u}(x,k)\ +\ \hat{v}(x,k) \\
& = & -\i k \int_{z\in E(r,\theta)}
f(z)\Bigl[ k^2\int_{\R^n}q(y)\,\Phi_k(y,z)\,\Phi_k(x,y)\,dy\Bigr. \\
& & \qquad \qquad \Bigl. +\ (q(z)+1)\,\Phi_k(x,z) \Bigr] ds(z) \\
& = & R\bigl[(f(\cdot)\,L(x,\cdot,k))\bigr](r,\theta) \,,
\end{eqnarray*}
where $R[f](r,\theta)$ is the $(n-1)$-dimensional Radon transform of $f$ in direction $(r,\theta)$ and
\begin{equation} \label{eq:B}
L(x,z,k)\ =\ (-\i k)^3\int_{\R^n}q(y)\,\Phi_k(y,z)\,\Phi_k(y,x)\,dy\ -\
\i k(q(z)+1)\,\Phi_k(x,z)
\end{equation}
for $x\in\Gamma$, $z\in B$, $k\in\R$, $x\not=z$.
\smalf
In other words, we have
\begin{equation} \label{eq:Hs_total}
\boxed{\F^{-1}\bigl(R^{-1}(\hat{u}+\hat{v})\bigr)(x,z,t)\ =\ f(z)\,\check{L}(x,z,t)\,,}
\end{equation}
for $x\in\Gamma$, $z\in B$, $k\in\R$, $x\not=z$.
\smalf
Therefore, from the knowledge of $\hat{m}^{r,\theta}(x,k)=\hat{u}(x,k)+
\hat{v}(x,k)$ for all $x\in\Gamma$, $k\in\R$, $r\geq 0$, and $\theta\in
S^{n-1}$ we can determine
$$ f(z)\,\check{L}(x,z,t)=\bigl(f_0(z)+f_1(z)\bigr)\,\check{L}(x,z,t) $$
for all $x\in\Gamma$, $z\in B$ with $x\not=z$, and $t\geq 0$.

\subsection{3D Domain}
\label{sec:3D}

For $n=3$ we recall that
\begin{equation*}
\Phi_k(x,y)\ =\ \frac{e^{\i k\abs{x-y}}}{4\pi\abs{x-y}}\,,\quad x\not= y\,.
\end{equation*}
Thus from (\ref{eq:B}) it follows that
$$ L(x,z,k)\ =\ (-\i k)^3 \frac{1}{16\pi^2}\int_{\R^3}q(y)\,
\frac{e^{\i k(\abs{y-z}+\abs{x-y})}}{\abs{x-y}\abs{y-z}}\,dy\ -\
\i k(q(z)+1)\,\frac{e^{\i k\abs{x-z}}}{4\pi\abs{x-z}}\,. $$
Taking the inverse Fourier transform with respect to $k$ gives
\begin{eqnarray*}
\check L(x,z,t) & = & \frac{1}{32\pi^3}\,\partial_{ttt} \left(
\int_{\R^3}\frac{q(y)}{\abs{x-y}\abs{y-z}}\intr
e^{\i k(\abs{y-z}+\abs{x-y}-t)}\,dk\,dy \right) \nonumber \\
& & +\ (q(z)+1)\,\frac{1}{8\pi^2\abs{x-z}}\,\partial_t \left(\intr
e^{\i k(\abs{x-z}-t)}\,dk \right) \nonumber \\
& = & \frac{1}{32\pi^3}\,\partial_{ttt} \left(\int_{\R^3} \delta(\abs{y-z}+\abs{y-x}-t)
\frac{q(y)}{\abs{z-y}\abs{y-x}} ds(y) \right) \nonumber\\
& &  -\frac{q(z)+1}{8\pi^2\abs{x-z}}\,\delta'(\abs{x-z}-t)\;.
\nonumber
\end{eqnarray*}
\medlf
Twice Integration with respect to $t$ gives
\begin{eqnarray*}
\lefteqn{\int_0^{s_2} \int_0^{s_1} \check L(x,z,\tau)\,d\tau ds_1} \\
& = &
\frac{q(z)+1}{8\pi^2\abs{x-z}}H(s_2-\abs{x-z}) \\
& & +\ \partial_{s_2} \left(\underbrace{\int_{\R^3} \delta(\abs{y-z}+\abs{y-x}-s_2)
\frac{q(y)}{32\pi^3\abs{z-y}\abs{y-x}}\,ds(y)}_{=:{\mathcal N}(x,z,s_2)} \right)\,,
\end{eqnarray*}
where we used the Heavyside function $H(\tau)=1$ for $\tau>0$ and
$H(\tau)=0$ for $\tau<0$. The initial conditions for the second term vanish
since for $t=0$ the domain of integration (w.r.t. $y$) is empty.
\smalf
Integrating once more, we get
\begin{eqnarray}
\Psi(x,z,t) & := & \int_0^t\int_0^{s_2}\int_0^{s_1}\check{L}(x,z,\tau)\,d\tau\,
ds_1\,ds_2 \label{eq:psi-def} \\
& = & \frac{q(z)+1}{8\pi^2\abs{x-z}}\int_0^tH(s_2-\abs{x-z})\,ds_2 + {\cal N}[q](x,z,t)\nonumber\,.
\nonumber
\end{eqnarray}
Now, we note that
\begin{equation*}
\int_0^t H(s_2-\abs{x-z})\,ds_2\ =\ (t-\abs{x-z})^+\ :=\ (t-\abs{x-z})\,
H(t-\abs{x-z})\,,
\end{equation*}
and thus
\begin{equation} \label{eq:psi}
\Psi(x,z,t) = \frac{q(z)+1}{8\pi^2\abs{x-z}}(t-\abs{x-z})^+\ + {\cal N}[q](x,z,t)\,.
\end{equation}
Both terms vanish for $t<\abs{x-z}$.
\smalf
We recall that $M(x,z,t):=f(z)\Psi(x,z,t)$ is known from measurement data for
all $x\in\Gamma$, $z\in B$, and $t>0$. First we take $z\in B\setminus\Omega$.
Then $q(z)=0$ and $f(z)=f_0(z)$, thus $\Psi(x,z,t)=M(x,z,t)/f_0(z)$ is known and
$$ \frac{M(x,z,t)}{f_0(z)}\ -\ \frac{(t-\abs{x-z})^+}{8\pi^2\abs{x-z}}\ = {\cal N}[q](x,z,t)\,. $$
Again, the left hand side is known. Now we let $z$ tend to $x$. The limit on
the right hand side exists, thus also on the left hand side, and therefore for
all $x\in\Gamma$ and $t>0$ we have
\begin{eqnarray}
\lim\limits_{z\to x}\left[\frac{M(x,z,t)}{f_0(z)}-\frac{t}{8\pi^2\abs{x-z}}
\right] & = & \frac{1}{8\pi^3\,t^2}\int_{\R^3} \delta(\abs{y-x}-t/2)
q(y)\,ds(y) \nonumber \\
& = & \frac{1}{8\pi}\,{\mathcal M}_2[q](x,t/2)\,, \label{q2D_reconstruct}
\end{eqnarray}
where ${\mathcal M}_2[q]$ is the spherical mean operator (see Appendix).
\smalf
Thus the reconstruction algorithm is as follows:
{\tt
\begin{enumerate}
\item Calculate the product $f(z)\check{L}(x,z,t)$ from (\ref{eq:Hs_total}) for
all $x\in\Gamma$, \\ $z\in B$ and $t>0$. By integration (see
(\ref{eq:psi-def})), this yields the knowledge of
$M(x,z,t):=f(z)\Psi(x,z,t)$ for all $x\in\Gamma$, $z\in B$ and $t>0$.
\item Solve (\ref{q2D_reconstruct}) for $q$ in $\Omega$ by inverting the
spherical mean o\-per\-a\-tor.
\item Compute $\Psi(x,z,t)$ for all $x\in\Gamma$, $z\in\Omega$ and $t>0$ from
(\ref{eq:psi}).
\item Finally, compute $f$ from $f(z)=M(x,z,t)/\Psi(x,z,t)$ for $z\in\Omega$.
\end{enumerate}}
\begin{remark}
The operator ${\mathcal N}[q](x,z,t)$ is the \emph{rotational ellipsoidal mean
operator} with focal points $x$ and $z$. Thus the integral equation
(\ref{eq:psi}); that is,
\begin{equation*}
\Psi(x,z,t)\ =\ \frac{q(z)+1}{8\pi^2\abs{x-z}}(t-\abs{x-z})^+\ +\
\mathcal{N}[q](x,z,t)\,,
\end{equation*}
can, for instance, be considered as a fixed point equation for $q$ involving the
ellipsoidal mean operator and can be solved by the fixed point iteration
\begin{equation*}
\Psi(x,z,t)\ =\ \frac{q_n(z)+1}{8\pi^2\abs{x-z}}(t-\abs{x-z})^+\ +\
\mathcal{N}[q_{n-1}](x,z,t)\,, n=1,2,\ldots
\end{equation*}
Ellipsoidal mean operators have been studied in John's book \cite{Joh04}.
\end{remark}

\subsection{2D Domain}
\label{sec:2D}

In $\R^2$ we recall that the fundamental solution is given by
\begin{equation*}
\Phi_k(x,y)\ =\ \frac{\i}{4}\,H_0^{(1)}(k\abs{x-y})\,,\quad x\not= y\,.
\end{equation*}
We have to compute the inverse Fourier transform of the product
$\Phi_k(x,y)\,\Phi_k(z,y)$ (see (\ref{eq:B})) and use the convolution theorem.
First we have that
$$ \frac{i}{4}\,\F^{-1}\bigl(H_0^{(1)}(\cdot\abs{x-y})\bigr)(t)\ =\
\left\{\begin{array}{cl} \displaystyle
\frac{1}{2\pi\sqrt{t^2-\abs{x-y}^2}}\,, & t>\abs{x-y}\,, \\
0\,, & t<\abs{x-y}\,. \end{array}\right. $$
The convolution for the inverse Fourier transform is given by
the formula \\ $\F^{-1}(fg)=\check{f}\ast\check{g}$; that is,
$$ \F^{-1}(fg)(t)\ =\ \bigl(\check{f}\ast\check{g}\bigr)(t)\ =\
\int_{-\infty}^\infty\check{f}(\tau)\,\check{g}(t-\tau)\,d\tau\,. $$
Therefore,
$$ \left(\frac{i}{4}\right)^2\F^{-1}\bigl(H_0^{(1)}(\cdot\abs{x-y})
H_0^{(1)}(\cdot\abs{z-y})\bigr)(t) $$
$$ =\
\left\{\begin{array}{cl}
\frac{1}{(2\pi)^2}\int_{\abs{x-y}}^{t-\abs{y-z}}\!\!\!
\frac{1}{\sqrt{\tau^2-\abs{x-y}^2}}\,\frac{1}{\sqrt{(t-\tau)^2-\abs{z-y}^2}}\,
d\tau\,, &  t-\abs{y-z}>\abs{x-y}\,, \\[4mm]
0\,, &  t-\abs{y-z}<\abs{x-y}\,, \end{array}\right. $$
and thus from (\ref{eq:B})
\begin{eqnarray}
\lefteqn{\check{L}(x,z,t)} \label{eq:B1} \\
& = & \frac{1}{4\pi^2}\,\partial_{ttt} \left(
\int_{\R^2}q(y)\int_{\abs{x-y}}^{t-\abs{y-z}}\!\!\!
\frac{1}{\sqrt{\tau^2-\abs{x-y}^2}}\,\frac{1}{\sqrt{(t-\tau)^2-\abs{z-y}^2}}
\,d\tau\,dy\right) \nonumber \\[2mm]
& & +\ \frac{1}{2\pi}\,\partial_t\left(\frac{q(z)+1}{\sqrt{t^2-\abs{x-z}^2}}
H(t-\abs{x-z})\right)\,, \nonumber
\end{eqnarray}
for $t-\abs{y-z}>\abs{x-y}$ where $H$ denotes again the Heaviside function.
Again, as in 3D, both terms vanish for $t<\abs{x-z}$.
\medlf
Letting $z$ tend to $x$ (thus $t>\abs{x-z}$ for $z$ sufficiently close to $x$)
yields
\begin{eqnarray*}
\lefteqn{\check{L}(x,x,t)} \\
& = & \frac{1}{4\pi^2}\,\partial_{ttt} \left(
\int\limits_{2|y-x|<t}q(y)\int\limits_{\abs{x-y}}^{t-\abs{y-z}}\!\!\!
\frac{1}{\sqrt{\tau^2-\abs{x-y}^2}}\,\frac{1}{\sqrt{(t-\tau)^2-\abs{x-y}^2}}\,
d\tau\,dy \right) \\[2mm]
& & \quad +\ \frac{1}{2\pi}\,\partial_t \left(\frac{1}{t} \right).
\end{eqnarray*}
Note that $q(x)$ vanishes for $x\in\Gamma$. Using polar coordinates and denoting
\begin{equation*}
k(t,r)\ :=\ \left\{\begin{array}{cl}\displaystyle
\int_r^{t-r}\frac{1}{\sqrt{\tau^2-r^2}}\,\frac{1}{\sqrt{(t-\tau)^2-r^2}}
\,d\tau\,, & t>2r>0\,, \\[5mm]
0\,, & 0<t<2r\,, \end{array} \right.
\end{equation*}
we therefore get
\begin{equation*}
\int_0^{t/2}k(t,r)\,r\int_0^{2\pi}q\left(x+r\binom{\cos\rho}{\sin\rho}\right)\,
d\rho\,dr\ =\ 2\pi\int_0^{t/2}\mathcal{M}_1[q](x,r)\,k(t,r)\,r\,dr
\end{equation*}
with the one dimensional circular mean operator $\mathcal{M}_1$. Using the
notations
\begin{eqnarray}
\tilde{k}(t,r) & := & r\,k(2t,r)\quad\mbox{and} \nonumber \\
\Phi(x,t) & := & 4\pi^2\int_0^t\int_0^{s_2}\int_0^{s_1}\check{L}(x,x,\tau)\,
d\tau\,ds_1\,ds_2\ -2\pi\bigl(t\log(t)-t\bigr)\,,\quad\mbox{ } \label{eq:Phi}
\end{eqnarray}
we get, for fixed $x$, the Volterra integral equation for
$\mathcal{M}_1[q](x,\cdot)$:
\begin{equation} \label{eq:volterra}
\Phi(x,t)\ =\ \int_0^t\tilde{k}(t,r)\,{\mathcal M}_1[q](x,r)\,dr\,,\quad
t\in(0,T]\,.
\end{equation}
We express $k$ as a complete elliptic integral: We introduce the variable
$\phi\in(-\pi/2,\pi/2]$ by $\tau=\frac{t}{2}+\left(\frac{t}{2}-r\right)
\sin\phi$. Then $d\tau=\left(\frac{t}{2}-r\right)\cos\phi\,d\phi$ and, with
$t'=t/2$,
\begin{eqnarray*}
\tau^2-r^2 & = & (\tau-r)(\tau+r) \\
           & = &\ \bigl[(t'-r)+(t'-r)\sin\phi\bigr]\,
\bigl[(t'+r)+(t'-r)\sin\phi\bigr]\,, \\
(t-\tau)^2-r^2 & = & (t-\tau-r)(t-\tau+r)\\
& = & \bigl[(t'-r)-(t'-r)\sin\phi\bigr]\,\bigl[(t'+r)-(t'-r)\sin\phi\bigr]\,.
\end{eqnarray*}
For the product we conclude that
\begin{eqnarray*}
\lefteqn{\bigl[\tau^2-r^2\bigr]\bigl[(t-\tau)^2-r^2\bigr]} \\
& = & \bigl[(t'-r)^2-(t'-r)^2\sin^2\phi\bigr]\,
\bigl[(t'+r)^2-(t'-r)^2\sin^2\phi\bigr] \\
& = & (t'-r)^2\cos^2\phi\,\bigl[(t'+r)^2-(t'-r)^2\sin^2\phi\bigr]\,.
\end{eqnarray*}
Therefore,
\begin{eqnarray*}
k(t,r) & = & \int_{-\pi/2}^{\pi/2}\frac{d\phi}
{\sqrt{(t'+r)^2-(t'-r)^2\sin^2\phi}} \\
& = & \frac{2}{t'+r}\int_0^{\pi/2}\frac{d\phi}
{\sqrt{1-\left(\frac{t'-r}{t'+r}\right)^2\sin^2\phi}}\,,
\end{eqnarray*}
and thus
$$ \tilde{k}(t,r)\ =\ \rho\,K(1-\rho)\quad\mbox{with}\quad
\rho=\frac{2r}{t+r}\,, $$
and where
$$ K(\alpha)\ =\ \int_0^{\pi/2}\frac{d\phi}
{\sqrt{1-\alpha^2\sin^2\phi}}\,,\quad|\alpha|<1\,, $$
denotes the complete elliptic integral.
\medlf
Let $\delta>0$ such that $\delta<\operatorname{dist}(\Gamma,\Omega)$. Since
${\mathcal M}_1[q](x,\tau)=0$ for $\tau\leq\delta$ we can consider the integral
equation (\ref{eq:volterra}) on $[\delta,T]$; that is,
\begin{equation} \label{eq:volterra-2}
\boxed{\Phi(x,t)\ =\ \int_\delta^t\tilde{k}(t,r)\,\mathcal{M}_1[q](x,r)\,dr\,,
\quad t\in[\delta,T]\,.}
\end{equation}
Now we note that $\tilde{k}(t,t)=K(0)=\pi/2$ and $\tilde{k}\in C^1(\Delta)$
where $\Delta=\{(t,r)\in[\delta,T]^2:r\leq t\}$. Therefore, we can transform
the equation of the first kind into one of the second kind by differentiating
(\ref{eq:volterra-2}). This yields
\begin{equation} \label{eq:volterra-3}
\boxed{
\Phi(x,t)\ =\ \frac{\pi}{2}\,{\mathcal M}_1[q](x,t)\ +\ \int_\delta^t
\frac{\partial}{\partial t}\tilde{k}(t,r)\,{\mathcal M}_1[q](x,r)\,dr\,,\quad
t\in[\delta,T]\,.}
\end{equation}
Volterra integral equations of the second kind with smooth kernel and non-vanishing diagonal have a unique
solution, see, e.g., \cite[Section~3.3]{Kress1999}.
\smalf
After calculating the spherical mean operator ${\mathcal M}_1[q]$, it can be inverted by using standard inversion formulas, which exist for
various geometries: For $B = B^2_R(0) \subset \R^2$, analytical reconstruction formulas have been derived by Finch, Haltmeier, Rakesh
\cite{FinHalRak07}. For a general domain~$\Omega$, Kunyansky reduced in~\cite{Kun07b} the reconstruction problem to the determination of
the eigenvalues $\lambda_k$ and normalized eigenfunctions $u_k$, $\|u_k\|_2=1$, of the Dirichlet Laplacian $-\Delta$ on $\Omega$ with
zero boundary conditions. Recently Palamadov \cite{Pal11} presented a general approach leading to reconstruction algorithms for various geometries.
For more details see also the Appendix.
\smalf
Thus the reconstruction algorithm is as follows:
{\tt
\begin{enumerate}
\item Calculate the product $M(x,z,t):=f(z)\check{L}(x,z,t)$ from
(\ref{eq:Hs_total}) for all $x\in\Gamma$, $z\in B$, $t>0$. Compute $\Phi(x,t)$
from (\ref{eq:Phi}) for
all $x\in\Gamma$ and $t>0$.
\item Solve (\ref{eq:volterra-3}) for calculating $\mathcal{M}_1[q]$.
\item Invert the spherical mean operator.
\item From (\ref{eq:B1}) determine $\check{L}$ by differentiation.
\item Finally use $f(z)=M(x,z,t)/\check{L}(x,z,t)$ to reconstruct $f(z)$ for
$z\in\Omega$.
\end{enumerate}}

\section*{Acknowledgements}
The research in this paper was initialized at the Program on Inverse Problems and
Applications at MSRI, Berkeley, during the Fall of 2010. It was followed up at the Special Semester at the Newton Institute 
in 2011.

\section{Appendix}
\subsection*{Notation}
$B_R^n(x)$ denotes the ball of radius $R$ and center $x$ in $\R^n$.

\subsection*{Bessel Functions}

The Bessel functions, Neumann functions and Hankel functions of order zero are
defined as
\begin{eqnarray*}
J_0(z) & = & \frac{1}{2\pi}\int_{-\pi}^\pi e^{\i z\sin(\tau)}\,d\tau\quad
\text{for }z\in\cmplx\,,\qquad\mbox{and} \\
Y_0(z) & = & \frac{2}{\pi}\,\bigl[\gamma+\ln(z/2)\bigr]\,J_0(z)\ -\
\frac{2}{\pi}\sum_{\ell=0}^\infty\frac{a_{\ell+1}}{(\ell!)^2}\,
\left(\frac{z}{2}\right)^{2\ell}\quad\mbox{for $z\in\cmplx$ with }
\Im\geq 0\,, \\
H_0^{(1)} & = & J_0\ +\ \i\,Y_0\,,
\end{eqnarray*}
where $\gamma$ denotes Euler's constant, $a_1=-\gamma$, and
$a_n=-\gamma+\sum_{\ell=1}^{n-1}\ell^{-1}$ for $n\geq 2$.
\smalf
From this we conclude that $H^{(1)}_0(-x)=\overline{H^{(1)}_0(x)}$ for $x\in\R$.
\subsection*{The Fourier Transform}
In this paper we use the following notation:
\begin{equation*}
\hat{f}(k)\ :=\ \F f(k)\ :=\ \intr f(t)\,e^{\i k t}\,dt\,,
\end{equation*}
denotes the Fourier transform of $f$ and
\begin{equation*}
\check{f}(t)\ :=\ \F^{-1}f(t)\ :=\ \frac{1}{2\pi}\intr f(k)\,e^{-\i t k}\,dk\,,
\end{equation*}
the inverse Fourier transform. For $f\in L^1(\R)$ or $f$ from the Schwarz
space $\mathcal{S}$ they are defined in the classical sense, for $f\in L^2(\R)$
they are defined by extension using the Parseval formula, and for tempered
distributions they are defined by duality.
\medlf
\textbf{Examples:}
\begin{eqnarray*}
\hat\delta(k)& = & \F\delta(k)\ =\ \intr e^{\i kt}\,\delta(t)\,dt\ \equiv\
e^{\i k0}\ =\ 1\,, \\
\check\delta(t) & = & \F^{-1}\delta(t)\ =\ \frac{1}{2\pi}\intr e^{-\i kt}\,
\delta(k)\,dk\ \equiv\ \frac{1}{2\pi}\,.
\end{eqnarray*}
This shows that
\begin{equation*}
\hat{1}(k)\ =\ 2\pi\delta(k)\quad\mbox{and}\quad
\check{1}(t)\ =\ \frac{1}{2\pi}\intr e^{-\i kt}\,dk\ =\ \frac{1}{2\pi}\,
\delta(t)\,.
\end{equation*}
Furthermore, we have (see \cite{Oberh1990})
$$ \frac{i}{4}\,H_0^{(1)}(ka)\ =\ \int\limits_a^\infty\frac{\exp(ikt)}
{2\pi\sqrt{t^2-a^2}}\,dt\quad\mbox{for }k>0\,. $$
From $H^{(1)}_0(-x)=\overline{H^{(1)}_0(x)}$ for all $x\in\R$ we conclude that
this formula holds for all $k\in\R$ with $k\not=0$. The right hand side is a
Fourier transform. Therefore,
$$ \mathcal{F}^{-1}\left(\frac{i}{4}\,H_0^{(1)}(a\cdot)\right)\ =\
\left\{\begin{array}{cl}
\displaystyle\frac{1}{2\pi\sqrt{t^2-a^2}}\,, & t>a\,, \\ 0\,, & t<a\,.
\end{array}\right. $$

\subsection*{Radon Transform}
In $\R^n$, $E(r,\theta)$ denotes the $(n-1)-$ dimensional hyperplane with orientation $\theta \in S^{n-1}$ and distance $r$ from origin.
\smalf
Let $f: \R^n \to R$ with support in $\Omega$, then the $(n-1)$-dimensional Radon transform of $f$ is defined by
\begin{equation*}
R_{r,\theta}(f)\ =\ R[f](r,\theta)\ =\ \int_{E(r,\theta)} f(x) ds(x)\;.
\end{equation*}
Thus $R[f]$ is a function from $\R \times S^{n-1}$ into $\R$.

\subsection*{Spherical Mean Operator}

In $\R^n$ the \emph{spherical mean operator} is defined as follows (see e.g. \cite{Eva98}):
\begin{equation*}
{\mathcal M}_{n-1}[u](x,r) := \frac{1}{\abs{S^{n-1}}} \int_{S^{n-1}} u(x+ry)\,ds(y) \text{ for } x \in \R^n, r > 0\;.
\end{equation*}
The spherical mean operator can be written as
\begin{equation*}
{\mathcal M}_{n-1}[u](x,r) = \frac{1}{r^{n-1}\abs{S^{n-1}}} \int_{\R^n} \delta(\abs{x-z}-r) u(z)\,dz \text{ for } x \in \R^n, r > 0\;.
\end{equation*}
In particular
\begin{equation*}
\begin{aligned}
{\mathcal M}_{2}[u](x,r) &= \frac{1}{4\pi^2}\int_{S^2} u(x+ry)\,ds(y) =
\frac{1}{4\pi^2 r^2} \int_{\R^3} \delta(\abs{x-z}-r) u(z)\,dz\,,\\
{\mathcal M}_{1}[u](x,r) &= \frac{1}{2\pi}\int_{S^1} u(x+ry)\,ds(y) =
\frac{1}{2\pi r} \int_{\R^2} \delta(\abs{x-z}-r) u(z)\,dz\;.
\end{aligned}
\end{equation*}
There exist a variety of reconstruction formulas for averages over circular and spherical means:
\subsubsection*{in 2D}
\begin{itemize}
\item
For averaged data on the sphere of Radius $R$, that is on  $\partial \Omega = \partial B^2_R(0) \subset \R^2$, analytical reconstruction
formulas have been derived by Finch, Haltmeier, Rakesh \cite{FinHalRak07} and read as follows
\begin{equation*}
f(\xi) = \frac{1}{2 \pi} \Delta_{\xi} \left(\int_{S^1} \int_0^{2R} r (\mathcal M_1[f])(\theta,r) \log \abs{r^2 -\abs{ \xi - \theta}^2 } \,d r \,ds(\theta)\right)
\end{equation*}
and
\begin{equation*}
f(\xi) = \frac{1}{2\pi} \int_{S^1} \int_0^{2}\left( \partial_r r \partial_r \mathcal M_1[f] \right)(\theta,r)
\log \abs{r^2 - \abs{\xi - \theta}^2} d r d s(\theta).
\end{equation*}
\item
For a general domain~$\Omega$, Kunyansky reduced in~\cite{Kun07b} the reconstruction problem to the determination of the eigenvalues $\lambda_k$ and normalized eigenfunctions $u_k$, $\|u_k\|_2=1$, of the Dirichlet Laplacian $-\Delta$ on $\Omega$ with zero boundary conditions:
\begin{align*}
\Delta u_k(\xi)+\lambda_ku_k(\xi)&=0,\quad\xi\in\Omega, \\
u_k(\xi)&=0,\quad\xi\in\partial\Omega.
\end{align*}
Indeed, if $(\xi,\eta)\mapsto G_{\lambda_k}(|\xi-\eta|)$ is a free-space rotationally invariant Green's function of the Helmholtz equation
and $n(\xi)$ denotes the outer unit normal vector of $\partial\Omega$ at $\xi\in\partial\Omega$, then
\begin{equation*}
f(\xi) = 2\pi\sum_{k=0}^\infty\tilde M_ku_k(\xi),
\end{equation*}
where
\[ \tilde M_k = \int_{\partial\Omega}\int_0^\infty r\mathcal M_1[f](\eta,r)G_{\lambda_k}(r)\left<\nabla u_k(\eta),n(\eta)\right>d rd s(\eta). \]
\end{itemize}
\subsubsection*{in 3D}
We refer to the survey of Finch and Rakesh \cite{FinRak09} (see also \cite{FinRak07}). In this article three reconstruction formulas are documented, which are
\begin{equation*}
\begin{aligned}
f(x) &= - \frac{1}{2 \pi R_0} \int_{\partial B_{R_0}^3} \frac{\partial_{tt} \left. \left(t^2 (\mathcal M_2[f])(x_0,t)\right) \right|_{t=\abs{x-x_0}}}{\abs{x-x_0}}\,ds(x_0)\,,\\
f(x) &= - \frac{1}{2 \pi R_0} \int_{\partial B_{R_0}^3} \frac{\left. \left( \partial_{t} t \partial_t t (\mathcal M_2[f])\right)(x_0,t) \right|_{t=\abs{x-x_0}}}{\abs{x-x_0}}\,ds(x_0)\,,\\
f(x) &= - \frac{1}{2 \pi R_0} \Delta \left(\int_{\partial B_{R_0}^3} \abs{x-x_0} \mathcal M_2[f])(x_0,\abs{x-x_0})\,ds(x_0)\right)\;.\\
\end{aligned}
\end{equation*}
Here $x \in B_{R_0}^3$ and $f$ is compactly supported in this set.

\def\cprime{$'$} \providecommand{\noopsort}[1]{}

%

\end{document}